\documentclass[10pt]{article}

\usepackage{amssymb}
\usepackage{latexsym}
\usepackage{amsmath}
\usepackage{amsfonts}
\usepackage{graphicx}
\usepackage{times}
\usepackage{subfig}
\usepackage{float}

\newtheorem{definition}{\bf Definition}[section]

\newtheorem{Thm}[definition]{\bf Theorem}
\newtheorem{Prop}[definition]{\bf Proposition}
\newtheorem{Cor}[definition]{\bf Corollary}

\def\R#1{{\operatorname{R}(#1)}}
\def\S#1{{\operatorname{S}(#1)}}
\def\T#1{{\operatorname{T}(#1)}}
\def\Q#1{{\operatorname{Q}(#1)}}

\title{Bounds on the differential of a graph}

\author{Ludwin A. Basilio$^{1}$, Jes\'us Lea\~nos$^{1}$, Omar Rosario$^{2}$ and Jos\'e M. Sigarreta$^{2}$
	\\
	\\
	$^1${\small Academic Unit of Mathematics, Autonomous University of Zacatecas. Paseo la
		Bufa, int. Calzada Solidaridad, 98060 Zacatecas, Zacatecas, Mexico}\\
	{\small ludwin.ali@gmail.com}
	\\
	$^2${\small Faculty of Mathematics, Autonomous University of Guerrero.
		Carlos E. Adame 5, Col. La Garita, Acapulco, Guerrero, Mexico}\\
	{\small orosario@uagro.mx}
	\\
	$^2${\small  Faculty of Mathematics, Autonomous University of Guerrero.
		Carlos E. Adame 5, Col. La Garita, Acapulco, Guerrero, Mexico}
	\\
	{\small josemariasigarretaalmira@hotmail.com}}

\begin{document}
	
	\begin{center}
		{\LARGE The differential on Graph Operator $\R{G}$}\\[12pt]
		{\large Ludwin A. Basilio$^{1}$,  Jes\'us Lea\~nos$^{1}$, Omar Rosario$^{2}$, \\[5pt] and Jos\'e M. Sigarreta$^{2}$}\\[12pt]
		
		$^1${\small Academic Unit of Mathematics, Autonomous University of Zacatecas. Paseo la Bufa, int. Calzada Solidaridad, 98060 Zacatecas, Zacatecas, Mexico} \\

		$^2${\small Faculty of Mathematics, Autonomous University of Guerrero.
			Carlos E. Adame 5, Col. La Garita, 39650, Acapulco, Guerrero, Mexico}\\[4pt]
	
		{\small e-mail: ludwin.ali@gmail.com; jesus.leanos@gmail.com; omarrosarioc@gmail.com; josemariasigarretaalmira\@@hotmail.com}\\[3pt]
			
	\end{center}
\begin{abstract}
Let $G=(V(G), E(G))$ be a simple graph with vertex set $V(G)$ and edge set
$E(G)$. Let $S$ be a subset of $V(G)$, and let $B(S)$ be the set of
neighbours of $S$ in $V(G) \setminus S$. The  differential
$\partial(S)$ of $S$ is the number $|B(S)|-|S|$. The maximum value
of $\partial(S)$ taken over all subsets $S\subseteq V(G)$ is the 
differential $\partial(G)$ of $G$. The graph $\R{G}$ is defined as the graph obtained from $G$ by adding a new vertex $v_e$  for each 
$e\in E(G)$, and by joining $v_e$ to the end vertices of $e$. 
In this paper we study the relationship between $\partial(G)$ 
and $\partial(R(G))$, and give tight asymptotic bounds for $\partial(R(G))$. 
We also exhibit some relationships between certain vertex sets of $G$ and $R(G)$ which involve
well known graph theoretical parameters. 
\end{abstract}

\vspace{0.1cm} \small{ {\it Keywords:}  Differential of a graph; Graph operators;
 Dominating set.
  
}

\section{Introduction}
Social networks, such as Facebook or Twitter, have served as a powerful
tool for communication and information disseminating. As
a consequence of their massive popularity, social networks now have a wide variety
of applications in the viral marketing of products and in political
campaigns. Motivated by their numerous applications, some
authors \cite{KeKlTa1,KeKlTa2} have proposed several influential maximization problems,
which share a fundamental algorith\-mic problem for information diffusion in social networks: the problem of determining the best
 group of nodes to influence the rest. As it was showed in \cite{BeFe}, the study of the differen\-tial $\partial(G)$ of a
graph $G$ could be motivated by such scenarios. Before moving on any further, let us introduce the basic definitions and notation 
that will be used in this paper. 

Throughout this paper, $G=(V(G),E(G))$ is a simple graph of order $n\geq 3$ with
vertex set $V(G)$ and edge set $E(G)$. For the rest of this section, let us assume that $v\in V(G)$ and that $S\subseteq V(G)$. 
The set of all vertices of $G$ that are adjacent to $v\in V(G)$ is the {\em neighborhood of} $v$ and is denoted by $N_G(v)$.  The set $N_G(v)\cup \{v\}$ is the 
{\em closed neighborhood of} $v$ and is denoted by $N_G[v]$. We define $N_G(S):=\bigcup_{v\in S}N_G(v)$ and $N_G[S]:=N_G(S)\cup S$.
We shall use $\overline{S}$ to denote $V(G)\setminus S$, and $N_{\overline{S}}(v)$ to denote the set of all vertices of $\overline{S}$ that are adjacent to $v$.  An \emph{external private neighbour of} $v\in S$
{\em with respect to} $S$ is a vertex $w\in N_{\overline{S}}(v)$ such
that $w\notin N_G(u)$ for any $u\in S\setminus\{v\}$. The set of all
external private neighbors of $v$ with respect to $S$ is denoted by epn$_G[v,S]$.
The {\em degree} $\delta_G(v)$ of $v$ is the number $|N_G(v)|$. As usual, the {\em minimum} and {\em maximum degree} 
of $G$ will be denoted by $\delta(G)$ and $\Delta(G)$, respectively. The subgraph of $G$ induced by $S$ will be denoted by $\langle S \rangle_G$. 
 
 From now on, we shall use $B_G(S)$ to denote $N_G(S)\setminus S$, and $C_G(S)$ to denote $V(G)\setminus (B_G(S)\cup S)$. Then $\{S,B_G(S),C_G(S)\}$ is a partition of $V(G)$. The {\em differential} $\partial_G(S)$ of $S\subseteq V(G)$ is the number $|B_G(S)|-|S|$, and the {\em differential} $\partial(G)$ of $G$ is the maximum  of $\{\partial_G(S): S\subseteq V(G) \}$. We call $S \subseteq V(G)$ a $\partial-set$ or a {\em differential set} if $\partial_G(S)=\partial(G)$.  If $S$ is a differential set, and it has minimum (maximum) cardinality  among all differential sets of $G$, then we call it a \emph{minimum (maximum) differential set}.  
From the definition of $\partial(G)$ it is easy to see that if $G$ is disconnected, and $G_1,\ldots ,G_k$ are its connected components, then
$\partial(G)=\partial(G_1)+\cdot\cdot\cdot+\partial(G_k)$. In view of this, all graphs under consideration in this paper will be connected  unless otherwise stated.

As long as no confusion can arise, we will omit the reference to $G$ in $N_G(v), N_G[v], \-N_G(S), \-N_G[S],$ 
epn$_G[v,S],  \delta_G(v), \langle S \rangle_G, B_G(S), C_G(S),\partial_G(S),$ etc.

The study of $\partial(G)$, together with a variety of other kinds of differentials of a set, started in 2006 \cite{MaHaHeHeSl}. 
Since then, the differential of a graph has received considerable attention 
\cite{ArCa,BaBeSi,BaBeLeSi,BaCaLeSi,BeDeMaSi,BeFe,BeFeSi,BeRoSi,CaCa,He,RoYo,Si, Si2}. 
The differential of a set $S\subseteq V(G)$ was previously  considered in \cite{GoHe}, where it was denoted by $\eta(S)$. 
In 2012  S. Bermudo and H. Fernau \cite{BeFeSi} showed a close relationship between $\partial(G)$ and the Roman domination number $\gamma_R(G)$ of $G$, namely,  $\partial(G)+\gamma_R(G)=n$. On the other hand, the minimum differential of an independent set was considered in \cite{Zh}, and the 
 $B$-differential (or enclaveless number) of $G$, defined as $\psi(G):= \max\{ |B(S)| : S \subseteq V(G)\}$, was  investigated in \cite{MaHaHeHeSl,Sl1}. 
 
A {\em graph operator} is a mapping $F : {\Lambda} \rightarrow {\Gamma}$, where $\Lambda$ and  $\Gamma$ are families of graphs. The study of how a graph invariant changes when a graph operator is applied is a classical pro\-blem in graph theo\-ry. 
In 1943, J. Krausz \cite{Kr} introduced four graph operators denoted by  
$\Q{G}, \R{G}, \S{G},$ and $\T{G}$. As can be seen in \cite{Bi,HaNo,MeReRoSi,Pr,RaLo,YaYaYe}, these operators have received a good deal of attention from  both sides, theoretical and practical. Our interest in such graph operators derives from their numerous applications, and from the fact that the class of graphs 
$\Gamma^*$ that can be constructed from another class of graphs, say $\Lambda$, by successive applications of such graph operators is much larger than 
$\Lambda$. The ultimate goal of our project is to give sharp results on the differential of any element of $\Gamma^*$, and this paper is a first step in that direction.   

\begin{figure}[h]
	\centering
	\subfloat{\includegraphics[scale=0.16]{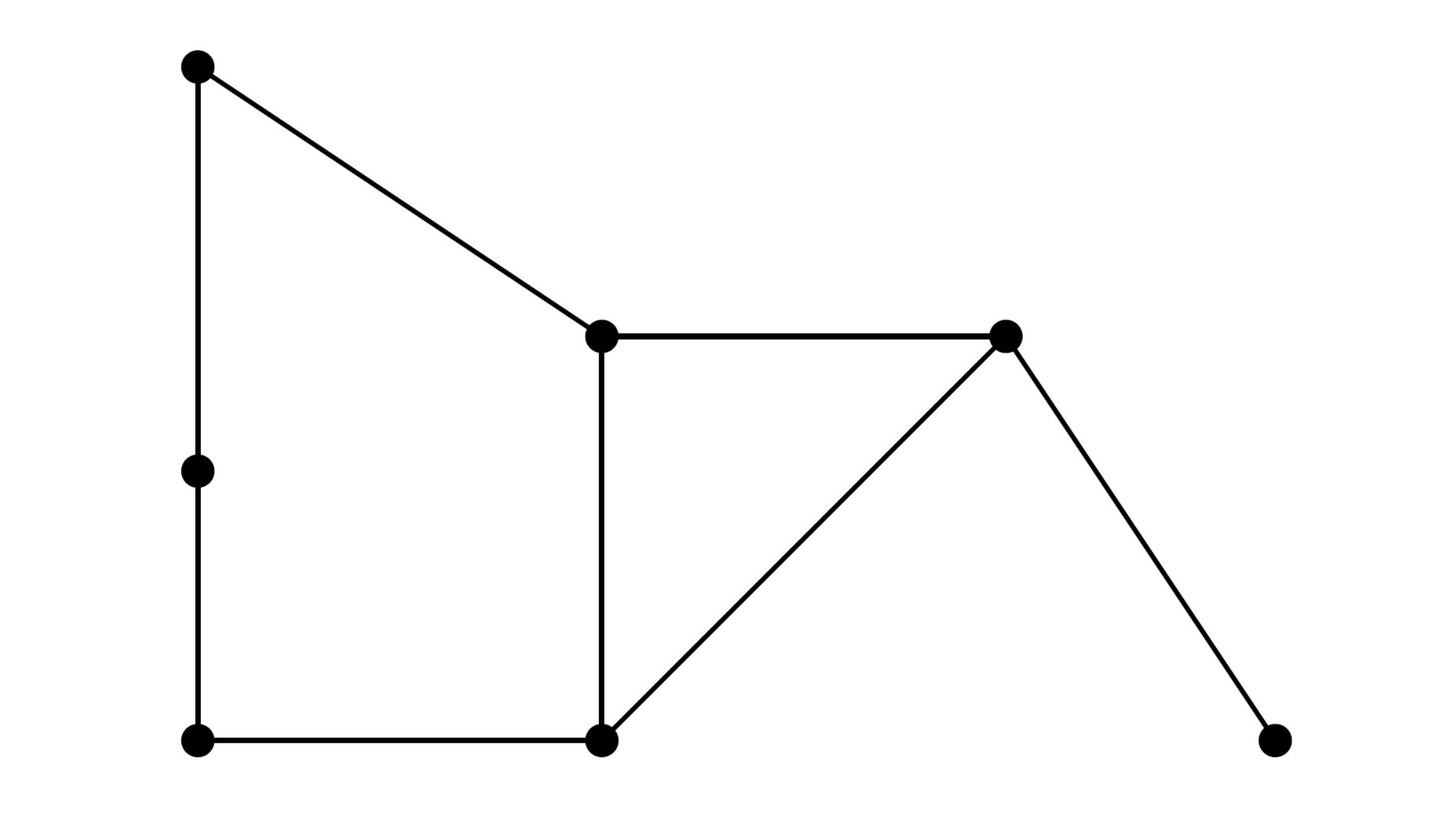}}
	\subfloat{\includegraphics[scale=0.26]{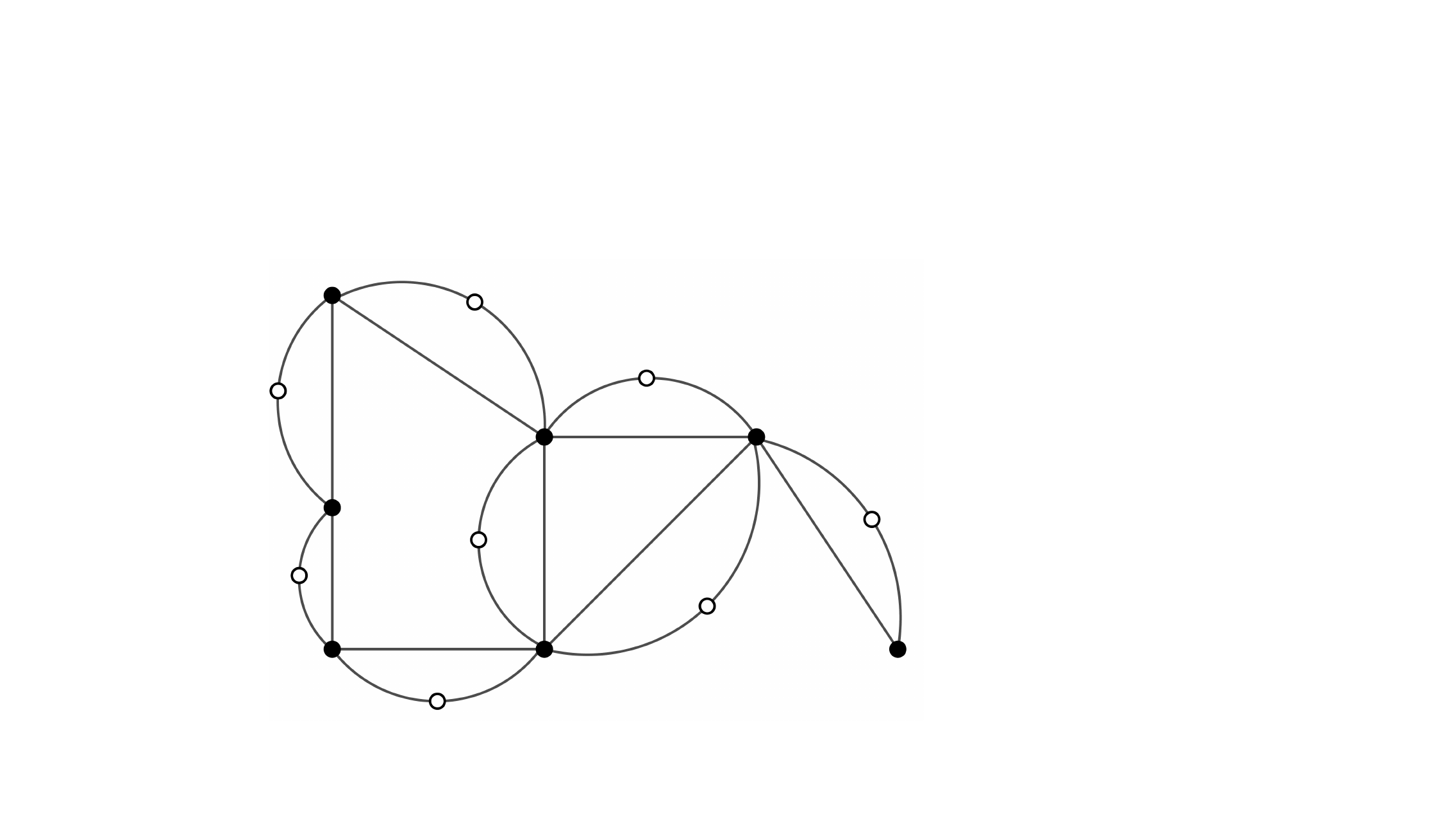}}
	\caption{The graph on the right is the corresponding $\R{G}$ of the graph $G$ on the left.}\label{fig:RG}
\end{figure}

In this work we deal only with the graph operator $\R{G}$. Following \cite{CvDoSa}, $\R{G}$ is defined as the graph obtained from $G$ by adding a new vertex 
$v_e$ for each edge $e$ of $G$, and by joining each $v_e$ to the end vertices of $e$. See Figure \ref{fig:RG} for an example. We begin our project with the study of $\R{G}$ because it shares structural properties with each of the other three operators. In view of this, we naturally hope that the present work  gives useful knowledge for the study of the rest of the operators.      

 \section{The differential of the graph $\R{G}$}
 
In this section we report some results on $\partial(\R{G})$. Most of such results show relations between certain vertex sets of $G$ and 
$\R{G}$. Each vertex set involved in these results is either a differential set, a dominating set, or a vertex cover set in at least one of $G$ or $\R{G}$.

 We start by listing some basic properties of $\R{G}$, which can be deduced easily from the definition of $\R{G}$. 
We use $E_n$ to denote the graph of order $n\geq 3$ with no edges. 

\begin{Prop}\label{obs-R(G)} 
Let $G=(V(G),E(G))$ be a graph of order $n$, and let $\R{G}$ be as above. Then
\begin{itemize}
\item[i)] $|V(\R{G})|=|V(G)|+|E(G)|.$  
\item[ii)] $|E(\R{G})|=3|E(G)|.$
\item[iii)] $G$ is an induced subgraph of $\R{G}.$
\item[iv)] $G \cong \R{G}$ if and only if $G \cong E_n$.
\item[v)] If $v\in V(G)$, then $\delta_{\R{G}}(v)=2\delta_G(v)$.  
\item[vi)] $G$ is connected and simple if and only if $\R{G}$ is connected and simple. 
\end{itemize} 
\end{Prop}

In view of Proposition \ref{obs-R(G)} $iii)$, we can partition $V(\R{G})$ into two subsets: the set of vertices $V$ formed by the vertices of $\R{G}$ that are in $G$, and $U:=V(\R{G})\setminus V$. Thus the subgraph $\langle V \rangle$ of $\R{G}$ is isomorphic to $G$, $|U|=|E(G)|$, and any vertex $u$ in $U$ has exactly two neighbours, say $v_1$ and $v_2$, such that $v_1,v_2\in V$ and $v_1v_2\in E(G)$. Such a partition $\{V,U\}$  will be called the {\em canonical partition}  of 
$V(\R{G})$. 
 
We recall that a {\em vertex cover} of $G$  is a subset $S\subseteq V(G)$ such that every edge of $G$ has at least one end vertex in $S$.  The 
{\em vertex conver number} of $G$ is the size of any smallest vertex cover in $G$ and is denoted by $\tau(G)$. A subset  $S\subseteq V(G)$ is a {\em dominating set} of $G$ if each vertex in $V(G)$ is in $S$ or is adjacent to a vertex in $S$. The {\em domination number} of $G$ is the minimum cardinality of a dominating set of $G$ and is denoted by $\gamma(G)$. For more information on domination in graphs see \cite{T1,T2}. Similarly, a subset  $S\subseteq V(G)$ is an {\em independent set} of $G$ if each any two distinct vertices in $S$ are not adjacent in $G$. The {\em independence number} of $G$ is the maximum cardinality of an independent set of $G$ and is denoted by $\alpha(G)$.




\begin{Prop}\label{prop:gammaV}
Let $\{V,U\}$ be the canonical partition of $V(\R{G})$. Then $V$ contains a minimum dominating set of $\R{G}$. 
\end{Prop}

\noindent \textit{Proof.} Let $W$ be a dominating set of $\R{G}$ such that $\gamma(\R{G})=|W|$. Since if $W \subseteq V$ then there is nothing to prove, we assume that  
$W\cap U$ has at least one vertex $u$. We know that $u$ has exactly two neighbours $v_1$ and $v_2$ such that $v_1,v_2\in V$ and $v_1v_2\in E(G)$. We can assume that $\{v_1,v_2\}\cap W=\emptyset$, as otherwise $W \setminus \{u\}$ is a dominating set of $\R{G}$, contradicting the minimality of $W$. Then $W':=(W \setminus \{ u \}) \cup \{ v_1\} $ is a minimum dominating set of $\R{G}$ with $|W'\cap U|=|W\cap U|-1$.  Clearly, continuing in this way, we can construct the required dominating set of $\R{G}$. 
\hfill \hfill $\square$

\begin{Prop}\label{prop:same-size}
Let $\{V,U\}$ be the canonical partition of $V(\R{G})$. If $D$ is a differential set of $\R{G}$, then $V$ contains a differential set $S$ of $\R{G}$ such that $|S|=|D|.$
\end{Prop}

\noindent \textit{Proof.}
Let $D$ be a subset of $V(\R{G})$ such that $\partial(\R{G})=\partial(D)$. Since if $D \subseteq V$ then there is nothing to prove, we assume that  
$D\cap U$ has at least one vertex $u$. Again, we know that $u$ has exactly two neighbours $v_1$ and $v_2$  
such that $v_1,v_2\in V$ and $v_1v_2\in E(G)$. Then $N(u) \subseteq N(v_i)$ for any $i\in \{1,2 \}$. This and the fact that  $\partial(\R{G})=\partial(D)$ imply
that neither $v_1$ nor $v_2$ belongs to $D$. Let  $D'=(D\setminus \{u\})\cup \{v_1\}$. Then $|D'|=|D|$, $B(D)\subseteq B(D')$, and so 
$\partial(D')=\partial(\R{G})$. Also note that  $|D'\cap U|=|D\cap U|-1$. Continuing in this way, we can construct the required $S$. 
\hfill \hfill $\square$

\begin{Prop}\label{Prop:domdif}
Let $\{V,U\}$ be the canonical partition of $V(\R{G})$. If $S \subseteq V$ is a differential set of $\R{G}$, then $S$ can be extended to a differential set $S'$ of $\R{G}$ such that: (i) $S'\subseteq V$, and  (ii) $S'$ is a dominating set of $G$.
\end{Prop}

\noindent \textit{Proof.}
 Let $S$ be as in the hypothesis of Proposition \ref{Prop:domdif}. We recall that $\{S, B(S),C(S)\}$ is a partition of $V(\R{G})$. If $S$ is a dominating set of $G$, then take $S'=S$ and we are done. Then we can assume that there exists a vertex $v_1\in C(S)\cap V$. Since $G$ is a connected graph of order at least $2$, then $v_1$ has at least one neighbour $v_2\in V$. If $v_2\in C(S)$, then the unique vertex $u\in U$ such that $N(u)=\{v_1,v_2\}$ must belong to $C(S)$. But then we have, 
 $\partial(S\cup \{v_1\})>\partial(S)=\partial(\R{G})$, a contradiction. Then we can assume that  all the neighbours of $v_1$ in $V$ are in $B(S)$. In particular, we have that $v_2\in B(S)$, and then $u\in C(S)$. Note that  $\partial(S\cup \{v_1\})=\partial(S)=\partial(\R{G})$. Continuing in this way, we can construct the required $S'$. 
\hfill \hfill $\square$

\begin{Cor}\label{Cor:clean}
Let $\{V,U\}$ be the canonical partition of $V(\R{G})$. Then $V$ contains a diffe\-rential set $S$ of $\R{G}$ such that $S$ is a dominating set of $G$.
\end{Cor}

\begin{figure}[H]
\centering

	\subfloat{\includegraphics[width=7.5cm]{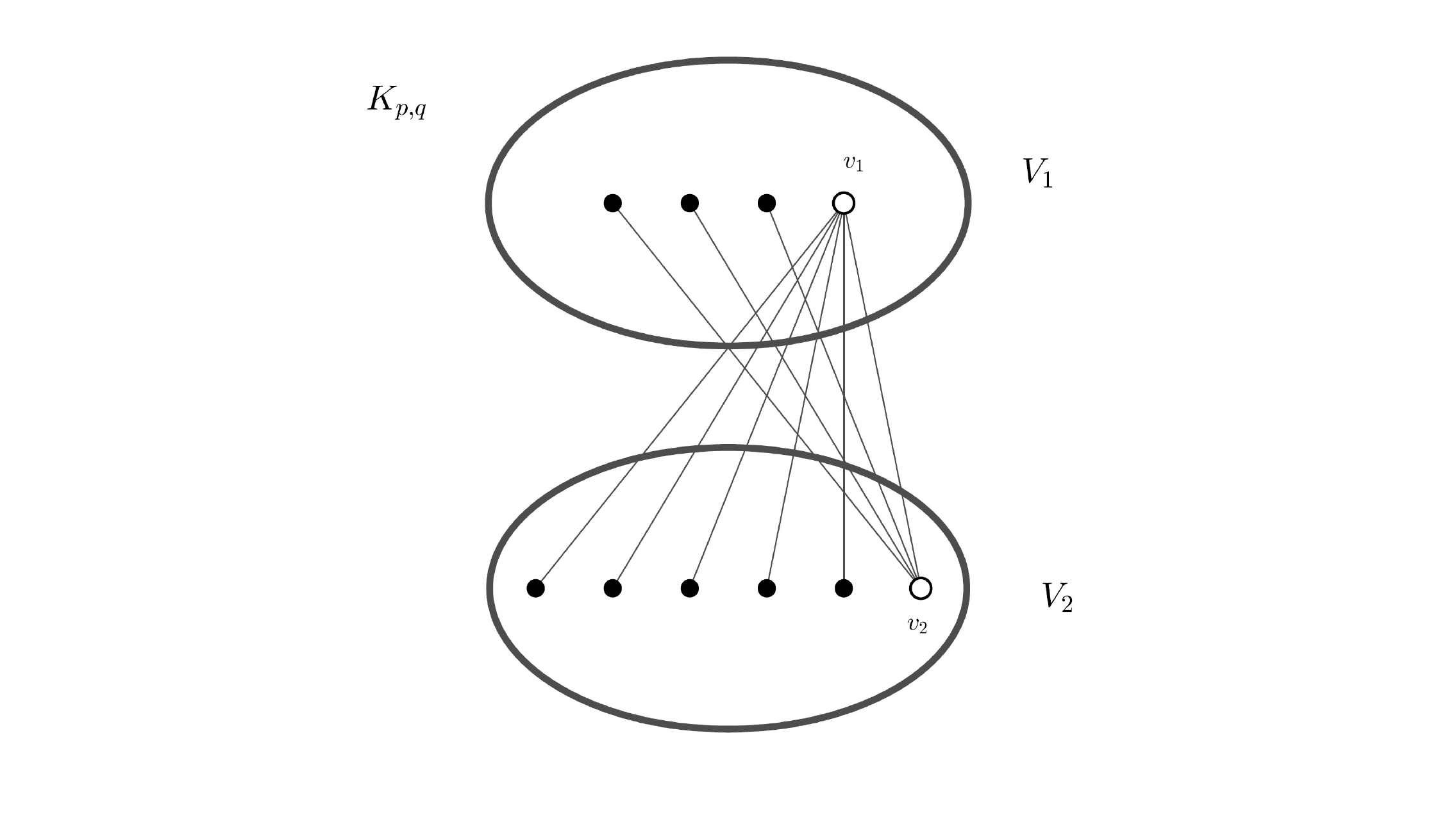}}
	\subfloat{\includegraphics[width=7.5cm]{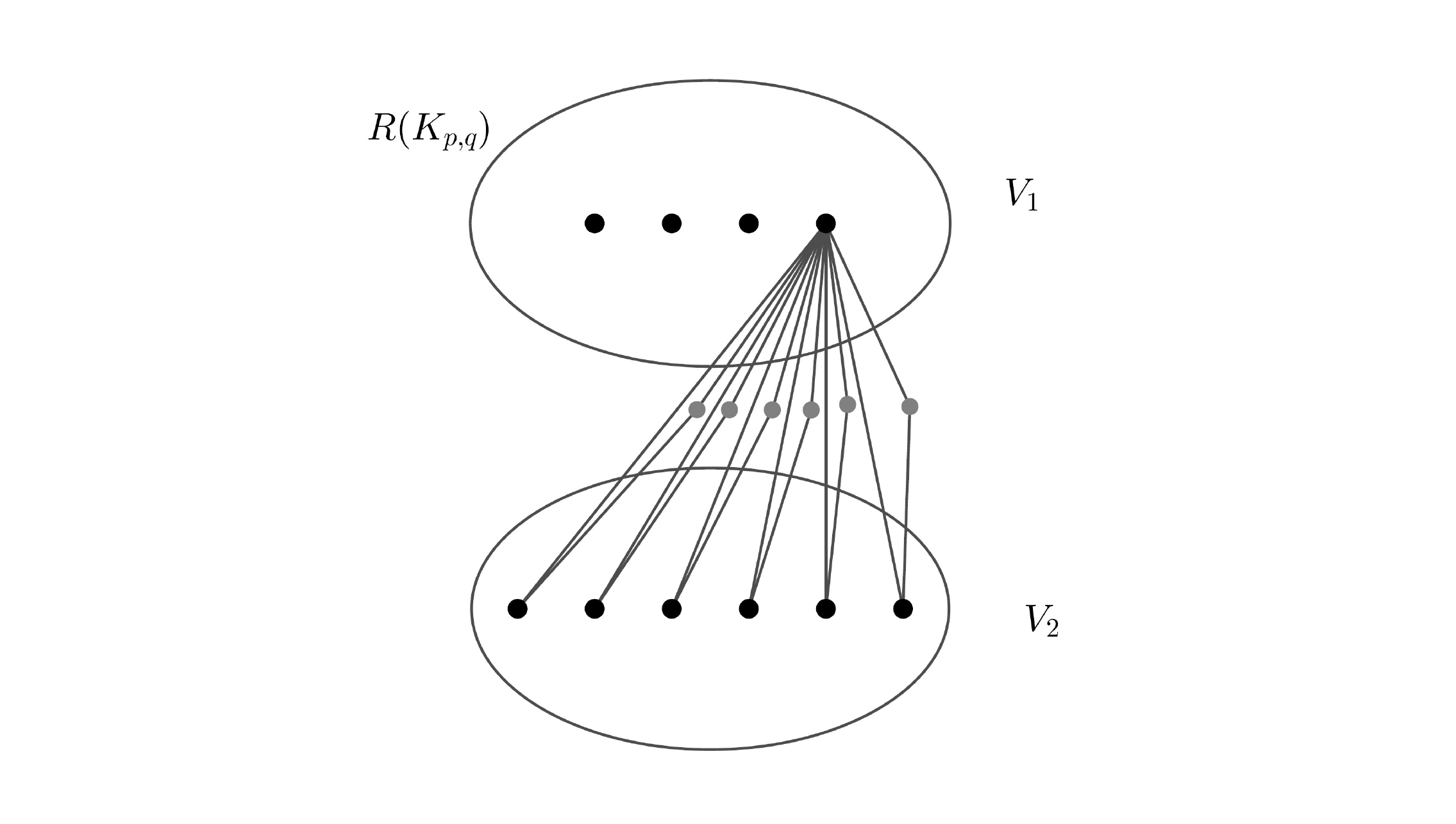}}
	\caption{The complete bipartite graph $G=K_{p,q}$ shows that the structure of a differential set of $G$ can be very different from the structure of a differential set of $\R{G}$.}\label{fig:bipartitacompleta}
\end{figure}

We remark that the set $S$ guaranteed by Corollary \ref{Cor:clean} is not necessarily a minimum domi\-nating set. Suppose that $S$ is the smaller part of the vertex bipartition of $V(K_{p,q})$. As we shall see later, if $|S|=p<q$, then such an $S$ is the only differential set of  $\R{K_{p,q}}$. Accor\-ding to Corollary \ref{Cor:clean}, $S$ is a dominating set for $K_{p,q}$, but it is far from being a minimum dominating set of $K_{p,q}$ (see Figure \ref{fig:bipartitacompleta}).
  
\begin{Prop}\label{prop:dominating}
Let $\{V,U\}$ be the canonical partition of $V(\R{G})$. If $S \subseteq V$ is a differential set of $\R{G}$ and $\delta(G)\geq 2$, then $S$ is a dominating set of $G$. 
\end{Prop}

\noindent \textit{Proof.} Let $S$ be as in the hypothesis of Proposition \ref{prop:dominating}. Since if $S$ is a dominating set of $G$ then there is nothing to prove, we assume that $V$ contains a vertex $w$ that is not adjacent to any vertex of $S$. Let $v_1$ and $v_2$ be two distinct neighbours of $w$ in $G$. Then $\R{G}$ has two distinct vertices $u_1,u_2 \in U$ such that $N(u_1)=\{w,v_1\}$ and $N(u_2)=\{w,v_2\}$. If $u_i \in N(S)$ for some $i\in \{ 1,2\}$, then we must have
 $v_i\in S$. This implies $w\in N(S)$, contradicting the choice of $w$.  Thus we can assume that $u_i \in C(S)$ for $i=1,2$, and so 
 $\partial(S \cup \{w\})>\partial(S)=\partial(\R{G})$. Since this inequality is false, we conclude that such a $w$ does not exist, as required. 
\hfill \hfill $\square$

\

\begin{figure}[h]
	\centering
	\includegraphics[width=10cm]{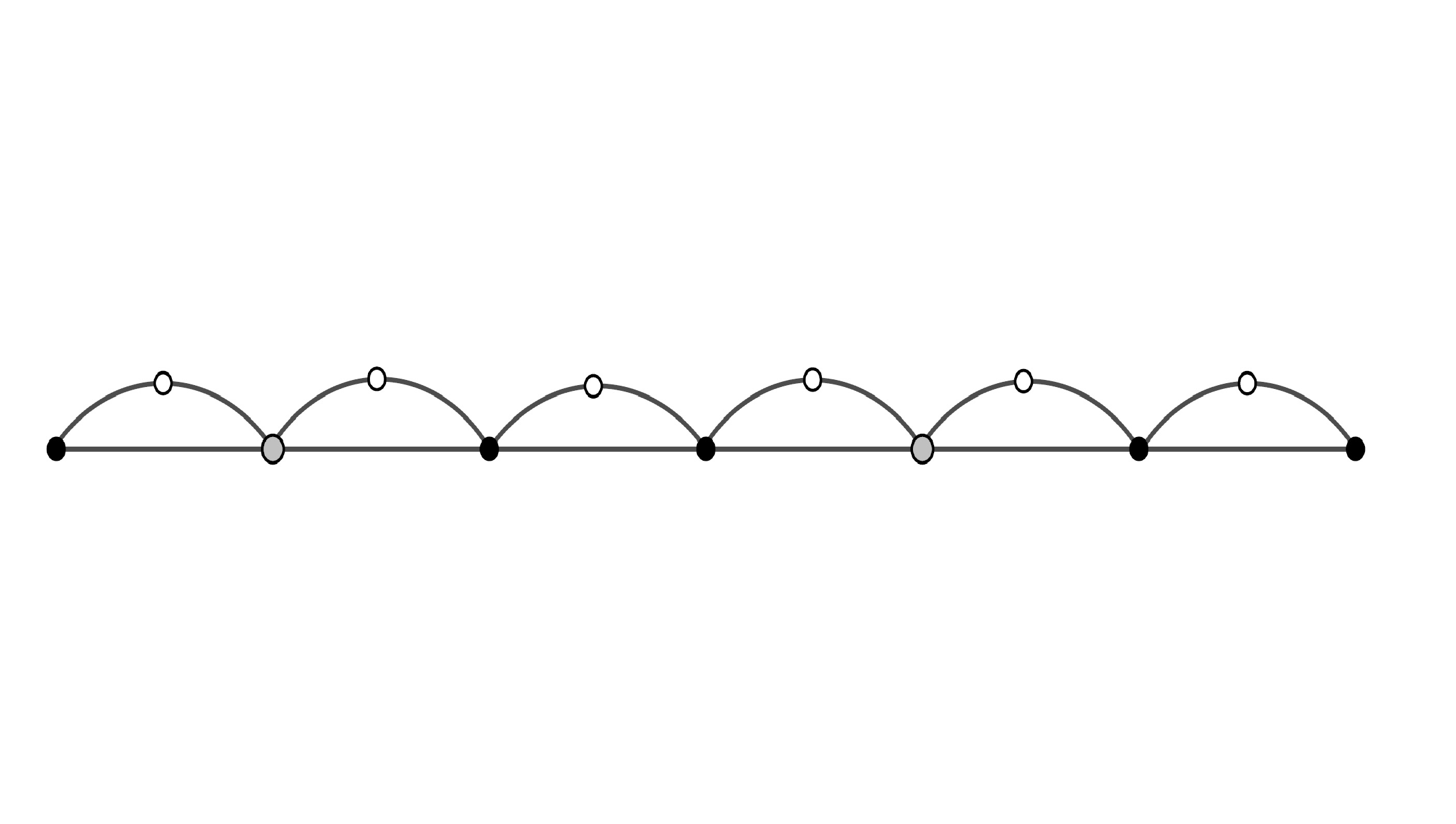}
	\caption{The minimum degree of $ P_7 $ is one, and the set $S$ of gray vertices form a differential set for both $ P_7 $ and $ \R {P_7} $. However, $S$ is not  a dominating set of $P_7$.}\label{fig:RGcamino}
\end{figure}

We remark that the reverse implication in Proposition \ref{prop:dominating} does not hold. Again, consider two vertices $x$ and $y$ in $K_{p,q}$ which belong to distinct parts of the vertex bipartition of $V(K_{p,q})$. Clearly, if $3\leq p\leq q$, then $K_{p,q}$ has minimum degree at least $2$ and $S=\{x,y\}$ is a (minimum) dominating set of $K_{p,q}$.  However,  $S$ is not a differential set of 
$\R{K_{p,q}}$. On the other hand, the set $S$ of gray vertices of the graph in Figure \ref{fig:RGcamino} shows that the condition $\delta(G)\geq 2$ in Proposition \ref{prop:dominating} is necessary.  

\begin{Cor}
Let $X \subseteq V(G)$ and $Y \subseteq V(\R{G})$ be differential sets of $G$ and $\R{G}$, respectively. If $\delta(G)\geq 2$, then
 $|Y|\geq |X|$.
\end{Cor}

Note that if $q>p\geq 3$, then any differential set of $K_{p,q}$ consists of exactly two vertices belonging to distinct parts of the vertex bipartition of $V(K_{p,q})$. On the other hand, we have proved in  Proposition \ref{cor:unique} that for the same range of values of $p$ and $q$, the only diffe\-rential set of $\R{K_{p,q}}$ is the smaller part of $V(K_{p,q})$. These facts show that the structure of a differential set of $G$ can be considerably distinct from the structure of a differential set of $\R{G}$. See Figure \ref{fig:bipartitacompleta} for an example.



\begin{Prop}{\cite{BeRoSi}}\label{prop:Dif,n-2}
	Let $G$ be a graph of order $n$ and maximum degree $\Delta$, then
	\begin{enumerate}
		\item[(a)]$\Delta = n-1$ if only if $\partial(G)=n-2$.
		\item[(b)]$\Delta = n-2$ if only if $\partial(G)=n-3$.
		\item[(c)]If $\Delta = n-3$, then $\partial(G)=n-4$.
	\end{enumerate}
\end{Prop}

\begin{Prop}
	Let $\{V,U\}$ be the canonical partition of $V(\R{G})$. If $m=|V(\R{G})|$, then  the following statements hold.  
		
	\begin{itemize}
		\item[(i)] $\partial(\R{G})=m-2$ if and only if $G=K_{1,n-1}$.
		
		\item[(ii)] $\partial(\R{G})=m-3$ if and only if $G=K_{1,n-1}\cup \{e\}$, where $e$ is an edge joining two leaf vertices of $K_{1,n-1}$. 
			\end{itemize}
\end{Prop}

\noindent \textit{Proof.} We recall  that in this paper $G$ is always a graph simple connected of order $n\geq 3$.  

We start by showing $(i)$. If $G =  K_{1,n-1}$,  then $\R{G}$ has a vertex $v$ such that $\delta(v)=m-1$, and by Proposition \ref{prop:Dif,n-2} $(a)$ we have that $\partial(\R{G})=m-2$.
	
	Now suppose that $\partial(\R{G})=m-2$. Then $\R{G}$ must have a vertex which is adjacent to any other vertex. This and the fact that  $n\geq 3$, imply that $v\in V$. In particular, $v$ must be adjacent to any other vertex of $G$, and so $G$ contains $K_{1,n-1}$ as subgraph. Since $v$ is adjacent to any other vertex of $\R{G}$, then  any edge of $G$ must belong to its subgraph $K_{1,n-1}$, and so $G=K_{1,n-1}$. 	
	
	Now we show $(ii)$. If $G=K_{1,n-1}\cup \{e\}$, where $e$ is an edge joining two leaf vertices of $K_{1,n-1}$, then the apex vertex $v$ of $K_{1,n-1}$
	has degree $m-2$, then by Proposition \ref{prop:Dif,n-2} $(b)$ and $(i)$ we have that $\partial(\R{G})=m-3$. 
	
	Now assume that  $\partial(\R{G})=m-3$. From $n\geq 3$ and the fact that $G$ is connected we have that $m\geq 5$. On the other hand, note that $\partial(\R{G})=m-3$ implies that the maximum degree of $\R{G}$ is $m-2$.
	
	 Let $v\in V(\R{G})$ with $\delta(v)=m-2\geq 3$. Since $\delta(u)=2$ for any $u\in U$, then we have that $v$ must belong to $V$. Since for each edge of $G$ with both ends in $V\setminus \{v\}$ there exists a unique vertex in $U$ which is not adjacent to $v$ and $\delta(v)=m-2$, then we have that 
	 $G$  has at most one edge with both ends in $V\setminus \{v\}$. On the other hand,  $\delta(v)=m-2$ implies that $\R{G}$ contains exactly one vertex, say $v_1$, which is not adjacent to $v$, and so  $C(\{v\})=\{v_1\}$. From the connectivity of $\R{G}$ we know that $ B(\{v\})$ contains a vertex $v_2$, 
	 which is adjacent to $v_1$. 
	 
	Note that if $v_2\in  U$, then we must have that $v_2$ is adjacent to $v$ in $\R{G}$, and hence that $v_1$ is adjacent to $v$ in $G$, contradicting that  $v_1\in C(\{v\})$. This implies that any neighbour of $v_1$ in $\R{G}$ must belong to $V\setminus \{v\}$, and hence that $v_1\in U$. Let $v_3$ be the other neighbour of $v_1$. Then the edge $e=v_2v_3$ belongs to $G$, and so  $G=K_{1,n-1}\cup \{e\}$, as required. 		 
\hfill \hfill $\square$

\begin{Prop}\label{cor:unique}
Let $p$ and $q$ be positive integers, such that $p+q\geq 4$ and $p<q$. Let $\{P,Q\}$ be the vertex bipartition of $K_{p,q}$, with $|P|=p$ and $|Q|=q$.  
Then $P$ is the only differential set of $\R{K_{p,q}}$. 
\end{Prop}

\noindent \textit{Proof.} 
Let $\{P, Q \}$ be the vertex partition of $V(K_{p,q})$ with  $|P|=p$ and $|Q|=q$. We recall that $p< q$. Since $P$ is a vertex cover of $K_{p,q}$, then
   $P$ is a dominating set of $\R{K_{p,q}}$. Then $\partial_{\R{K_{p,q}}}(P)=q(p+1)-p$, and so  $\partial(\R{K_{p,q}})\geq q(p+1)-p$. The assertion is easy to verify for $n=4,5,6$. Thus we shall assume that $n>6$.   
   
  Let $S$ be a differential set of $\R{K_{p,q}}$. By Proposition \ref{prop:same-size}, we can assume that $S\subseteq V$. Clearly, it is enough to show that $|S|\geq p$.
 Seeking a contradiction, suppose that $|S|\leq p-1$. Then each of $P$ and $Q$ has at least a vertex which is not in $S$. Let $v\in Q\setminus S$, and suppose that $P\setminus S$ has at least three vertices, say  $v_1, v_2,$ and $v_3$.  For $i=1,2,3$, let  $u_i$ be the only vertex of $U$ which is adjacent to  both $v$ and $v_i$. Note that the existence of such $u_1,u_2,u_3,v_1,v_2,v_3$ and $v$ implies, $\partial_{\R{K_{p,q}}}(S\cup\{v\})=\partial_{\R{K_{p,q}}}(S)+1$, a contradiction. Thus we conclude that $P$ has at most $2$ vertices which are not in $S$. By an analogous reasoning, we can conclude that $Q$ has at most $2$ vertices which are not in $S$. These imply $|S|\geq p+q-4$, and hence that $p-1\geq p+q-4$. But this last inequality implies, $3 \geq q \geq n/2$, a contradiction.

We know that if $S$ is a differential set of $\R{K_{p,q}}$, then $S\cap Q=\emptyset$. Thus we can assume that 
$S\subseteq P\cup U$.  We claim that  $S\cap U=\emptyset$. Suppose not, then there exists $u\in S\cap U$. Let $x\in P$ and $y\in Q$ be the only two neighbours of $u$ in $\R{K_{p,q}}$.
Thus $x$ cannot belong to $S$. 

On the other hand, since $q>p$ and $p+q\geq 4$, then $q\geq 3$. Let $y_1,y_2,\ldots ,y_q$ be the vertices of $Q$. For $i=1,2,\ldots, q$, let $u_i$ be the only vertex of $U$ such that
$N_{K_{p,q}}(u_i)=\{x,y_i\}$.  Since $x\notin S$ and $S\cap Q=\emptyset$, then $\{u_1,u_2,\ldots, u_q\}\cap B_{\R{K_{p,q}}}(S)=\emptyset$. Then $S':=S\setminus\{u\}\cup \{x\}$ satisfies the following:
$\partial_{\R{K_{p,q}}}(S')\geq \partial_{\R{K_{p,q}}}(S)+2$, which contradicts the choice of $S$. 
\hfill \hfill $\square$
 
 \
 
 We recall that the {\em wheel graph} $W_{n}$ of order $n\geq 4$ is formed by the cycle $C_{n-1}$, plus an additional vertex, say $v$, such that  $v$ is adjacent to any vertex in $V(C_{n-1})$.  
Usually, the vertex $v$ is called the {\em apex vertex} of $W_n$.

\begin{Prop}\label{prop:exactly} Let $p,q, n$ be positive integers, such that $p+q=n\geq 4$ and $p\leq q$. The following equalities hold: 
	\begin{itemize}
		\item[(i)] $\partial(\R{K_n})=\dfrac{n(n-1)}{2}-n+3$.
		
		\item[(ii)] $\partial(\R{W_n})=2n-3$.
		
		\item[(iii)] $\partial(\R{K_{p,q}})=q(p+1)-p$.
		
	\end{itemize}
\end{Prop}

\noindent \textit{Proof.} Let  $\{V,U\}$ be the canonical partition of $V(\R{G})$, where $G$ is $K_n, W_n,$ or $K_{p,q}$, depending on the case under consideration.

First we show $(i)$. Since $G=K_n$, then $|V|=n$ and $|U|=\frac{n(n-1)}{2}$. Let $S$ be a differential set of $\R{K_n}$. By Proposition \ref{prop:same-size}, we can assume that $S\subseteq V$. In particular, this implies that 
$1\leq |S|\leq n$.

We claim that $|S|\geq n-3$. Suppose not. Then $V\setminus S$ contains at least 4 vertices, say $v_1,v_2,v_3$ and $v_4$. Then these 4 vertices must belong to $B_{\R{K_n}}(S)$. For $i=2,3,4$, let
$u_{i}$ be the only vertex of $U$ which is adjacent to $v_1$ and $v_i$. Since $v_1,v_2,v_3,v_4\in B_{\R{K_n}}(S)$, then 
$u_2,u_3$ and $u_4$ are in $C_{\R{K_n}}(S)$. From these facts it follows that $\partial_{\R{K_n}}(S\cup \{v_1\})>\partial_{\R{K_n}}(S)$, contradicting that $S$ is a differential set of $\R{K_n}$.  Thus we must have that $|S|\in \{n-3, n-2,n-1,n\}$. 

On the other hand, the following is easy to check:  

\begin{equation}\label{eq:Kn}
\partial_{\R{K_n}}(S)= \left\{ 
\begin{array}{lcc}
             \frac{n(n-1)}{2}-n+3 &   if  & \hskip 1cm |S|=n-3,n-2,\\
             \\  \frac{n(n-1)}{2}-n+2 &  if & |S|=n-1,\\
               \\ \frac{n(n-1)}{2}-n  &  if & |S|=n.
                          \end{array}
   \right.
\end{equation}
From Equation \ref{eq:Kn} and the fact that these four are the only possibilities for such a set $S$, we can conclude that $\partial(\R{K_n})=\dfrac{n(n-1)}{2}-n+3$, as required.

Now we show $(ii)$.  Let $v$ be the apex vertex of $W_n$. Note that $|V(\R{W_n})|=3(n-1)+1=3n-2$ and also that $\partial_{\R{W_n}}(\{v\})=2(n-1)-1$. This last equality implies  
$\partial(\R{W_n})\geq 2n-3$. Since the case in which $n=4$ is easy to verify, we assume that $n\geq 5$. 

Let $S$ be a differential set of $\R{W_n}$, and let $\ell:=|S|>0$. By Proposition \ref{prop:same-size}, we can assume that $S\subseteq V$. Moreover, we claim that if $v\notin S$, then 
$S':=S\cup\{v\}$ is also a differential set of $\R{W_n}$. Indeed, suppose that  $v\notin S$, and that $\partial_{\R{W_n}}(S)>\partial_{\R{W_n}}(S')$. Then the last inequality implies that $V(C_{n-1})$ has at most one vertex not in $S$. This, the supposition that $v\notin S$, and $S\subseteq V$ imply that $|S|\geq n-2$. Hence $\partial_{\R{W_n}}(S)\leq ((3n-2)-1)-2|S|\leq 3n -3 - 2(n-2)=n+1$.  
Since $\partial(\R{S})\geq 2n-3$, then we have that $n+1\geq 2n-3$, or that $4\geq n$. Since $4\geq n$ contradicts $5\leq n$, we have that $S'$ is also a differential set of $\R{W_n}$. 

Let $v,v_{i_1},\ldots, v_{i_\ell}$ be the vertices of $S'$. We claim that $S'':=S'\setminus\{v_{i_\ell}\}$ is also a differential set of $\R{W_n}$.  Indeed, since epn$[v_{i_\ell},S']\leq 2$ and  
$|S''|=|S'|-1$, then we must have that $\partial_{\R{W_n}}(S')=\partial_{\R{W_n}}(S'')$, as required. Similarly, it can be deduced that $S''':=S''\setminus\{v_{i_{\ell-1}}\}$ is also a 
differential set of $\R{W_n}$.  Continuing in this way, we can conclude that $\{v\}$ is a differential set of $W_n$, and hence that $\partial(\R{W_n})=\partial_{\R{W_n}}(\{v\})=2n-3$, as required. 

By Proposition \ref{cor:unique}, we have $(iii)$.
\hfill \hfill $\square$

 

\

The following proposition shows a surprising relationships between the vertex cover number of $G$ and the domination number of $\R{G}$.

\begin{Prop}
Let $\{V,U\}$ be the canonical partition of $V(\R{G})$. Then $\tau(G)=\gamma(\R{G})$. 
\end{Prop}

\noindent \textit{Proof.} First we show that $\gamma(\R{G})\leq \tau(G)$. Let $X\subseteq V(G)$ be a vertex cover of $G$ such that $\tau(G)=|X|$. It suffices to show that $X$ is a dominating set of $\R{G}$.   Let $u$ be a vertex of $V(\R{G})$. Then we need show that $u\in X$ or that $\R{G}$ contains an edge $ux$ with $x\in X$. If $u\in X$ there is nothing to show. Then we can assume that $u\notin X$. 
If $u\in U$, then we know that $u$ has exactly two neighbours $v_1$ and $v_2$  
such that $v_1,v_2\in V$ and $v_1v_2\in E(G)$. Since $X$ is a vertex cover of $G$, then at least one of $v_1$ or $v_2$ must belong to $X$, and hence such a vertex is the required $x$. Thus we can assume that $u\in V$. Since $G$ is a connected graph of order $n\geq 3$, then $u$ has at least one neighbour in $V$, say $y$. From the fact that $X$ is a vertex cover of $G$, and the fact $u\notin X$ we have that $y$ must belong to $X$, and so $y$ is the required $x$.

Now we show that $\gamma(\R{G})\geq \tau(G)$. From Proposition \ref{prop:gammaV}, we know that $V$ contains a dominating set $S$ of $\R{G}$. Then it is suffices to show that $S$ is a vertex cover of $G$. Let $v_1v_2$ be an edge of $G$. From the definition of $\R{G}$ we know that $U$ contains a unique vertex $u$ such that $N(u)=\{ v_1,v_2\}$. Since such a vertex $u$ must be dominated by some vertex of $S$,
then at least one of $v_1$ or $v_2$ is in $S$, and hence $S$ is a vertex cover of $G$. This implies that $\tau(G)\leq \gamma(\R{G})$.
\hfill \hfill $\square$

\begin{Prop}\label{prop:vc+dif} 
Let $\{V,U\}$ be the canonical partition of $V(\R{G})$.  If $S\subseteq V$ is a  vertex cover of $G$ and
 $\partial(G)=\partial(S)$, then $S$ is a differential set of $\R{G}$.
\end{Prop}

\noindent \textit{Proof.} Let $S$ be as in the hypothesis of Proposition \ref{prop:vc+dif}. Since $S$ is a vertex cover of $G$, then $S$ is a dominating set for both $G$ and  $\R{G}$. As $S$ is a dominating set of $G$ and $\partial(G)=\partial(S)$, then $S$ must be a minimum dominating set of $G$. On the other hand, since $G$ is an induced subgraph of $\R{G}$, then $|S|=\gamma(G)\leq \gamma(\R{G})$. From the last inequality and the fact that $S$ is  a dominating set of $\R{G}$, we have that $S$ is also a minimum dominating set of $\R{G}$. Thus, $\gamma(\R{G})=|S|$.
 
Let $X \subseteq V(\R{G})$ such that  $\partial(\R{G})=\partial(X)$. Then $|X|\leq |S|$. By Corollary \ref{Cor:clean}, we can assume that 
$X$ is a dominating set of $G$ and that $X\subseteq V$. Since $S\subseteq V$ is a minimum dominating set of $G$, then $|X|\geq |S|$, and so $|X|=|S|=\gamma(\R{G})$. 

Let $m=|V(\R{G})|$. It is well known that if $Y \subseteq V(\R{G})$, then  $|B(Y)|\leq m-\gamma(\R{G})$. Then $\partial(\R{G})=\partial(X)=|B(X)|-|X|\leq m -\gamma(\R{G})-|X|\leq m- 2|S| =\partial(S) \leq \partial(\R{G})$.
\hfill \hfill $\square$

\begin{figure}[H]
	\centering
	\includegraphics[width=7.4cm]{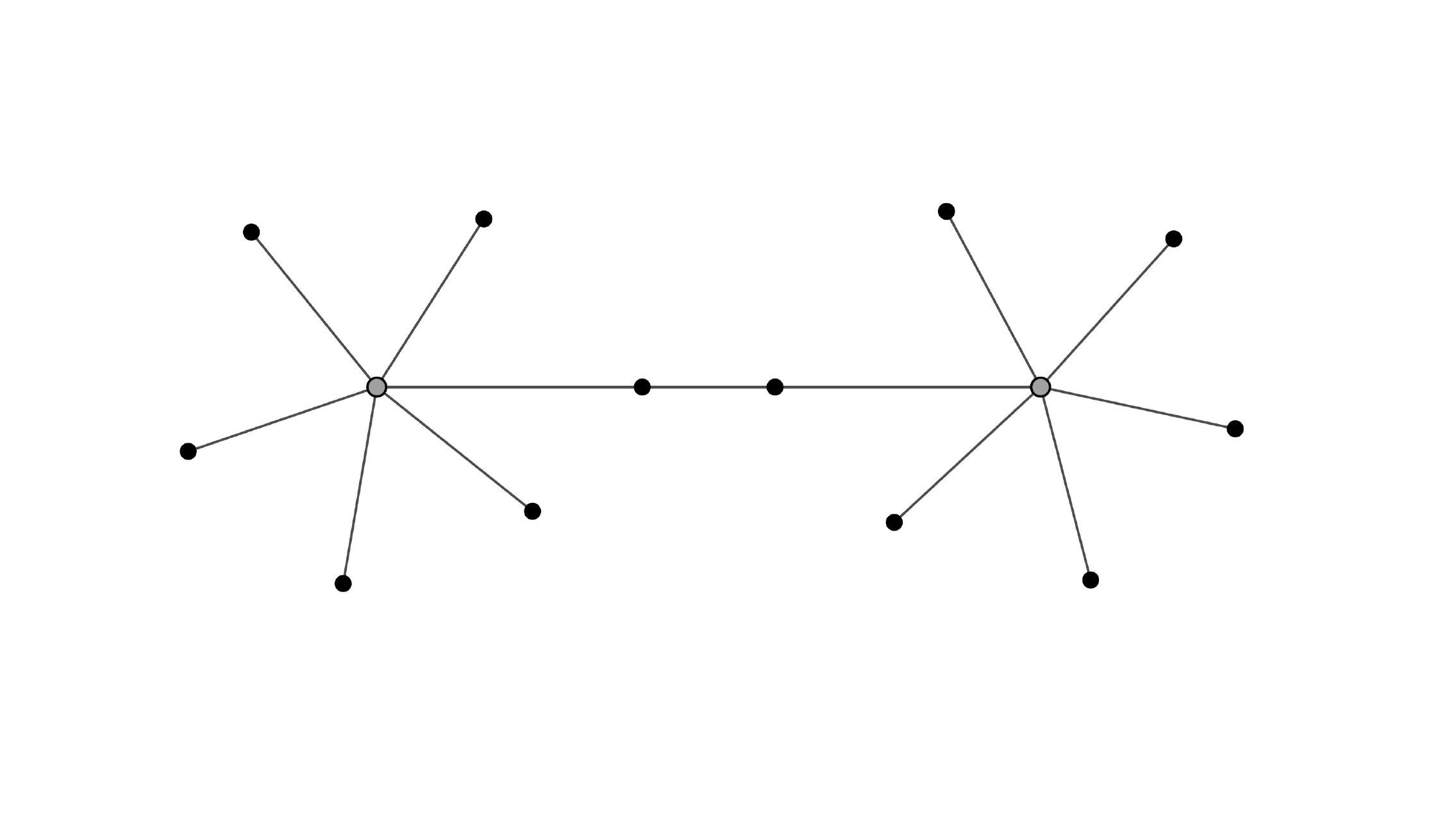}
	\includegraphics[width=7.4cm]{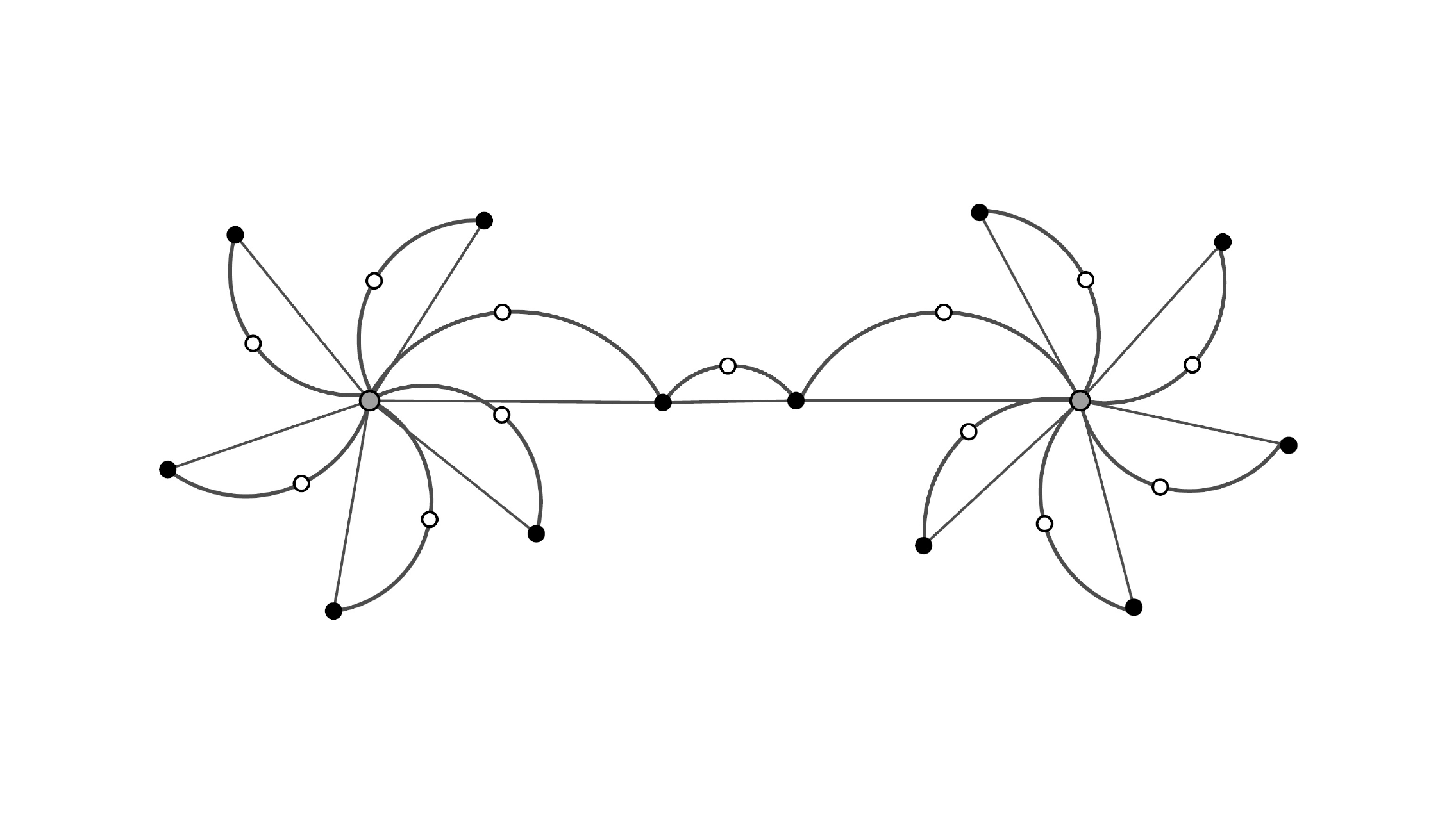}
	\caption{The gray vertices form a differential set of $G$ and $\R{G}$, but they are not a vertex cover of $G$.}\label{fig:dobestrella}
\end{figure}

\
Note that if $S$ is a differential set of $G$ and $\R{G}$, then $S$ is not necessarily a vertex cover of $G$.  See for instance the graph in  Figure \ref{fig:dobestrella}.

\section{General tight bounds  for $\partial(\R{G})$}
Our goal in this section is to show Theorem \ref{thm:main1}, which states general lower and upper bounds for $\partial(\R{G})$. 

We recall that a subset $S \subseteq V(G) $ is $k-dependent$ in $G$ if $\langle S \rangle$ has maximum degree at most $k$. For more about $k$-dependent sets we refer the reader to \cite{FaHeHeRa}.

\begin{Prop}\label{prop:leq2}
Let $\{V,U\}$ be the canonical partition of $V(\R{G})$. If $S\subseteq V$ is a diffe\-rential set of $\R{G}$, then $\langle B_{G}(S) \rangle$ is $2$-dependent.
Moreover, if $S$ is a maximal differential set of $\R{G}$, then $\langle B_{G}(S) \rangle$ is $1$-dependent.
\end{Prop}

\noindent \textit{Proof.}
We recall that $|V(G)|=n\geq 3$. Let $S\subseteq V$ be a differential set of $\R{G}$. Seeking a contradiction, suppose that $\langle B_{G}(S) \rangle$ has a vertex $v$ such that $\delta_{\langle B_{G}(S) \rangle}(v)\geq 3$. Let $v_1,v_2$ and $v_3$ be three distinct neighbours  of $v$ in $\langle B_{G}(S) \rangle$, and for $i=1,2,3$, let $u_i$ be the only vertex in $U$ that is adjacent to both $v$ and $v_i$. Then
$u_1,u_2,u_3 \in C_{\R{G}}(S) \cap N_{\R{G}}(v)$, and hence we have $\partial_{\R{G}}(S\cup \{ v\})\geq \partial_{\R{G}}(S)+1$, a contradiction.

Now suppose that, additionally, $S$ is maximal and that has $\langle B_{G}(S) \rangle$ a vertex $v$ of degree $2$. Let $v_1$ and $v_2$ be two distinct neighbours of 
$v$. For $i=1,2$, let $u_i$ be the only vertex in $U$ that is adjacent to both $v$ and $v_i$. Then
$u_1,u_2 \in C_{\R{G}}(S) \cap N_{\R{G}}(v)$, and hence we have $\partial_{\R{G}}(S\cup \{ v\})\geq \partial_{\R{G}}(S)$, which contradicts the maximality of $S$.
\hfill \hfill $\square$

If $S\subseteq V(G)$ is a differential set of $\R{G}$ of maximum cardinality, then we shall use $\mu(G)$ to denote  $|S|$. Similarly, we shall use $\lambda(G)$ to denote $|E(G)|-|V(G)|+2\alpha(G)$. Note that the connectivity of $G$ and $\alpha(G)\geq 1$ imply 
that $\lambda(G)\geq 1$.

\begin{Prop}\label{cor:leq2}
Let $\{V,U\}$ be the canonical partition of $V(\R{G})$. If $S\subseteq V$ is a maximum differential set of $\R{G}$, then $|C_{\R{G}}(S)|\leq \frac{|V(G)|-\mu(G)}{2}.$
\end{Prop}

\noindent \textit{Proof.}
Let $H$ be the  subgraph of $G$ induced by $B_{G}(S)$. We claim that $V\cap C_{\R{G}}(S)=\emptyset$. Suppose not, and let $v$ be a vertex in 
$V\cap C_{\R{G}}(S)$. Since $G$ is connected and $|V(G)|=n\geq 3$, then $v$ is adjacent to another vertex $v_1\in V$. 
Let $u_1$ be the only vertex of $U$ which is adjacent to both $v$ and $v_1$. Since $S$ is a differential set of $\R{G}$, then we must have that $v_1\in B_{\R{G}}(S)$ and that 
$u_1\in C_{\R{G}}(S)$. These last two conclusions imply that  $S\cup \{v\}$ is also a differential set of $\R{G}$, contradicting the maximality of $S$. Then, we can assume that 
$C_{\R{G}}(S)\subseteq U$. In particular, this implies that each $u\in C_{\R{G}}(S)$ corresponds to an edge of $H$, and hence  $|C_{\R{G}}(S)|$ is exactly the number of edges of $H$.  On the other hand, since $H$ has maximum degree at most $1$ by Proposition \ref{prop:leq2}, then $H$ has at most $\frac{|V(G)|-\mu(G)}{2}$ edges, as required.
\hfill \hfill $\square$

\begin{Thm}\label{thm:main1}
Let $\{V,U\}$ be the canonical partition of $V(\R{G})$. Then
$\lambda(G) \leq \partial(\R{G})\leq \lambda(G)+\lfloor\frac{|V(G)|-\mu(G)}{2}\rfloor$.  
\end{Thm}
\noindent \textit{Proof.} We recall that $|V(G)|=n\geq 3$. First we shall show $\lambda(G) \leq \partial(\R{G})$. Let $I$ be a maximum independent set of $G$, and let $S:=V\setminus I$. Then $\alpha(G)=|I|$ and $|S|=n-\alpha(G)$. 
Note that $B_{\R{G}}(S)=U\cup I$. Then 
$\partial(\R{G})\geq \partial_{\R{G}}(S)=|U|+|I|-|S|=|E(G)|+\alpha(G)-(n-\alpha(G))=|E(G)|+2\alpha(G)-n=\lambda(G)$, as required.

Now  we shall show that $\partial(\R{G})\leq \lambda(G)+\lfloor\frac{n-\mu(G)}{2}\rfloor$. Let $S$ be a maximum differential set of $\R{G}$. Then $\mu(G)=|S|$. Moreover, by Proposition \ref{prop:same-size} we can assume that $S\subseteq V$. Also note that the  maximality of $S$ implies that $C_{\R{G}}(S)\subseteq U$. 
Thus $V$ has exactly $n-\mu(G)$ vertices in $B_{\R{G}}(S)$. Let $F$ be the set of edges of $G$
which have both end vertices in  $V\setminus S$ and let $f:=|F|$. From Proposition \ref{prop:leq2} we know that $f\leq \lfloor\frac{n-\mu(G)}{2}\rfloor$, and that  $F$ is a matching of $G$. Let $S'$
be a set of vertices in $V$ which results by adding to $S$ exactly one vertex of each edge in $F$. Then $S'$ is a dominating set of $\R{G}$, and $I':=V\setminus S'$ is an independent set of $G$. Clearly, $\lambda(G)=|E(G)|+\alpha(G) -(n-\alpha(G))\geq |E(G)| + |I'| - (n-|I'|)=\partial_{\R{G}}(I')=\partial(\R{G})-f$.  Then $\lambda(G)\geq \partial(\R{G})-f\geq  \partial(\R{G})-\lfloor\frac{n-\mu(G)}{2}\rfloor$, as required.    
\hfill \hfill $\square$

For $r$ a positive integer, let us denote by $K'_{r,2r}$ the graph that results from the complete bipartite graph $K_{r,2r}$ when we add a matching of $r$ edges with both end vertices in the bigger part of the bipartition of $V(K_{r,2r})$.

Our next result shows that  both bounds in Theorem \ref{thm:main1} are tight. 

\begin{Prop}\label{prop:tight}
Let $r\geq 2$ be a positive integer. Then the following hold:
\begin{itemize}
\item[(i)] $\partial(\R{K_{r,2r}})= \lambda(K_{r,2r})$.

\item[(ii)] $\partial(\R{K'_{r,2r}})= \lambda(K'_{r,2r})+\left\lfloor\frac{3r-\mu(K'_{r,2r})}{2}\right\rfloor$.
\end{itemize}
\end{Prop}

\noindent \textit{Proof.} As usual, throughout this proof we assume that $\{V,U\}$ is the canonical partition of $V(\R{G})$, where $G$ is $K_{r,2r}$ or $K'_{r,2r}$ 
depending on the case under consideration.  

First we show $(i)$. Let $\{P,Q\}$ be the vertex bipartition of $V(K_{r,2r})$. We assume that $|P|=r$ and that $|Q|=2r$.  
Since $\alpha(K_{r,2r})=2r$, then $\lambda(K_{r,2r})=2r(r)-3r+4r=2r^2+r$. 
Since $P$ is vertex cover of $K_{r,2r}$, then $P$ is a dominating set of $\R{K_{r,2r}}$, and hence $\partial_{\R{K_{r,2r}}}(P)=2r(r+1)-r=2r^2+r$. In particular, we have that  $\partial(\R{K_{r,2r}})\geq \lambda(K_{r,2r})$.    
    
It remains to show that $P$ is a differential set of $\R{K_{r,2r}}$. Let $S$ be a differential set of $\R{K_{r,2r}}$. By Proposition \ref{prop:same-size}, we 
can assume that  $S\subseteq V$. 
Since $P$ is a dominating set of $\R{K_{r,2r}}$, then, in order to show that $P$ is a differential set of $\R{K_{r,2r}}$, it is enough to show that $|S|\geq r$.
 Seeking a contradiction, suppose that $|S|\leq r-1$. Then each of $P$ and $Q$ has at least a vertex which is not in $S$. In fact, $Q\setminus S$ has at leas $r+1\geq 3$ of such 
 vertices.  
 Let $v\in P\setminus S$, and let   $v_1, v_2,$ and $v_3$ be three vertices in $Q\setminus S$.  For $i=1,2,3$, let $u_i$ be the only vertex of $U$ which is adjacent to  both $v$ and $v_i$. Note that the existence of such $u_1,u_2,u_3,v_1,v_2,v_3$ and $v$ implies that $\partial_{\R{K_{r,2r}}}(S\cup\{v\})\geq \partial_{\R{K_{r,2r}}}(S)+1$, contradicting the choice of $S$. 
  
 Now we shall show $(ii)$. For brevity, we will use $G$ to denote $K'_{r,2r}$. Without loss of generality, we assume that 
\begin{itemize}
\item $P:=\{v_1,\ldots ,v_r\}, Q:=\{w_1,\ldots ,w_{2r}\}$, and
\item $E(G):=\{v_iw_j | v_i\in P \mbox{ and } w_j\in Q\} \cup \{w_iw_{i+r} | i=1,\ldots ,r\}$.
\end{itemize}
 Let $S$ be a maximum differential set of $\R{G}$. Then $\mu(G)=|S|$ and $\partial(G)=\partial_\R{G}(S)$. As before, by Proposition \ref{prop:same-size} we 
can assume that $S\subseteq V=P\cup Q$. It is a routine exercise to show that the domination number of $\R{G}$ is exactly $2r$. 
In view of this, we can assume that $|S|\leq 2r$. 
 
Each of the following equalities follows directly from the definitions involved: $(i)$ $\alpha(G)=r$, $(ii)$ $\lambda(G)=(2r(r)+r)-3r+2r=2r^2$, and $(iii)$ $\partial_{\R{G}}(P)= 2r(r)+|Q|-|P|=2r(r)+2r-r= 2r^2+r$. 

In particular, note that if $S=P$, then there is nothing to prove. Indeed, in such a case we have that $\mu(G)=r$, and hence $\partial_{\R{G}}(P)= 2r^2+r=2r^2+\lfloor\frac{3r-r}{2}\rfloor$, as required. 
Thus, from now on we assume that 
$S\neq P$. Let $P_1:=P\cap S, Q_1:=Q\cap S, P_0:=P\setminus P_1$, and $Q_0:=Q\cap S$. 

{\sc Claim 1.} $|Q_0|\geq 1$. Seeking a contradiction, suppose that $Q_0= \emptyset$. Then $Q\subseteq S$, and hence we must have that $Q=S$. This last implies that 
$\partial_{\R{G}}(S)=|U|+|P|-|S|=(2r(r)+r)+(r)-(2r)=2r^2<2r^2+r=\partial_{\R{G}}(P)$, which contradicts the choice of $S$.

{\sc Claim 2.} $|Q_1|\geq 1$. Seeking a contradiction, suppose that $Q_1=\emptyset$. Since $S\neq P$, then  $P_0\neq \emptyset$, and hence $|S|<|P|=r$. Then 
 $\partial_{\R{G}}(S)= (2r)|S|+|Q|-|S|=|S|(2r-1)+2r<r(2r-1)+2r=2r^2+r=\partial_{\R{G}}(P)$, contradicting the choice of $S$.

{\sc Claim 3.}  $|P_0|\geq 1$. Seeking a contradiction, suppose that $P_0= \emptyset$. This supposition and Claim 2 imply that $S=P\cup Q_1$ where $Q_1\neq \emptyset$. Let $U_1$ be the set of
vertices of $U$ which are adjacent to some vertex of $Q_1$ but not to a vertex of $P$. From the definition of $G$ and $\R{G}$ it is clear that $|Q_1|\geq |U_1|$.  Thus 
$\partial_{\R{G}}(S)=\left(2r(r)+|Q|-|Q_1|+|U_1|\right)-\left(|P|+|Q_1|\right)=(2r^2+|Q|-|P|)+(|U_1|-2|Q_1|)=\partial_{\R{G}}(P)+|U_1|-2|Q_1|.$ Since $|U_1|-2|Q_1|<0$, then
 $\partial_{\R{G}}(S)<\partial_{\R{G}}(P)$, which contradicts the choice of $S$. 
 
 {\sc Claim 4.}  $|P_1|\geq 1$. Seeking a contradiction, suppose that $P_1= \emptyset$, and hence that $P=P_0$. From Claim 1 we know that $Q_0$ contains at least one vertex, say $w$. 
  Let $U_w$ be the set of vertices of $U$ which have a neighbour in $P$ and the other in $\{w\}$. Since $|P|=r\geq 2$, then $|U_w|\geq 2$. Since no vertex in $U_w$ belongs to $B_{\R{G}}(S)$, then
 for $S':=S\cup \{w\}$ we have that $\partial_{\R{G}}(S')\geq \partial_{\R{G}}(S)$, contradicting the choice of $S$. 

 {\sc Claim 5.}  $|Q_0|=1$. Seeking a contradiction and considering Claim 1, we can suppose that $|Q_0|\geq 2$. From Claim 3 we know that $P_0$ contains at least one vertex, say $v$.
Let $U_v$ be the set of vertices of $U$ which have a neighbour in $\{v\}$ and the other in $Q_0$. Since $|Q_0|\geq 2$, then $|U_v|\geq 2$. Since no vertex in $U_v$ belongs to $B_{\R{G}}(S)$, then
 for $S':=S\cup \{v\}$ we have that $\partial_{\R{G}}(S')\geq \partial_{\R{G}}(S)$, contradicting the choice of $S$. 
 
 Claim 5 implies that $|Q_1|=2r-1$. This, together with $|S|\leq 2r$ and Claim 4 imply that $|P_1|=1$, and hence that $|S|=2r$.  Since $|Q_0|\geq 1$ and $|P_0|\geq 1$, then the set $U_{0,0}$ of vertices of $U$ which have a neighbour in $Q_0$ and the other in $P_0$ is nonempty. Moreover, since no vertex in $U_{0,0}$ belongs to $B_{\R{G}}(S)$, then $S$ is not a dominating set of $\R{G}$, but it  has cardinality $2r$. Thus $\partial_{\R{G}}(Q) > \partial_{\R{G}}(S)$, contradicting the choice of $S$.  
\hfill \hfill $\square$

\section{Acknowledgements} This research was partly supported by SEP(F-PROMEP-39/Rev-04)(Mexico). We also thank the referee for their valuable comments. 


\end{document}